\documentstyle
    [times,pstricks,pst-node,amssymb,amsthm,amsmath,twocolumn,aps,floats]
    {revtex}

\newcommand{\R}{{\mathbb{R}}}
\newcommand{\Z}{{\mathbb{Z}}}
\newcommand{\SO}{{\mathrm{SO}}}
\renewcommand{\O}{{\mathrm{O}}}
\newcommand{\ie}{\emph{i.e.}}
\newcommand{\Ie}{\emph{I.e.}}
\renewcommand{\tensor}{\otimes}
\renewcommand{\tilde}{\widetilde}
\newcommand{\fig}[1]{Figure~\ref{#1}}
\newtheorem{theorem}{Theorem}
\newtheorem{proposition}{Proposition}
\newtheorem{lemma}{Lemma}
\newtheorem{conjecture}{Conjecture}

\psset{linewidth=.5pt,dash=5pt 5pt}

\SpecialCoor

\begin{document}

\wideabs{
\title{Circumscribing constant-width bodies with polytopes}
\author{Greg Kuperberg}
\address{UC Davis \\ E-mail: \tt greg@math.ucdavis.edu}
\maketitle
\begin{abstract}
Makeev conjectured that every constant-width body is inscribed in the
dual difference body of a regular simplex.  We prove that homologically,
there are an odd number of such circumscribing bodies in dimension
3, and therefore geometrically there is at least one.  We show that
the homological answer is zero in higher dimensions, a result which
is inconclusive for the geometric question.  We also give a partial
generalization involving affine circumscription of strictly convex bodies.
\end{abstract}}

\begin{figure}[htb]
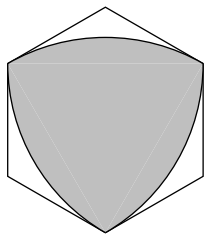

\begin{center}
\pspicture(-1.5,-1.5)(1.5,1.5)
\pspolygon(1.5;30)(1.5;90)(1.5;150)(1.5;210)(1.5;270)(1.5;330)
\pspolygon*[linecolor=lightgray](1.5;270)(1.5;30)(1.5;150)
\psarc[fillstyle=solid,fillcolor=lightgray](1.5;270){2.598}{60}{120}
\psarc[fillstyle=solid,fillcolor=lightgray](1.5;30){2.598}{180}{240}
\psarc[fillstyle=solid,fillcolor=lightgray](1.5;150){2.598}{300}{0}
\endpspicture
\end{center}
\caption{A Rouleaux triangle inscribed in a regular hexagon}
\label{f:rouleaux}
\end{figure}

Any set of diameter 2 in $\R^n$ is contained in a convex body of
constant width 2.  Consequently, if some polytope $P$ circumscribes every
convex body of constant width 2, it contains every set of diameter 2.
For example, every constant-width body in two dimensions is inscribed
in a regular hexagon (\fig{f:rouleaux}). A conjecture of Makeev
\cite{Makeev:polygons} generalizes this theorem to higher dimensions:

\begin{conjecture}[V. V. Makeev] Every constant width body in $\R^n$
is inscribed in a polytope similar to $D_n$, the dual of the difference
body of a regular simplex. \label{c:main}
\end{conjecture}

\begin{figure}[htb]
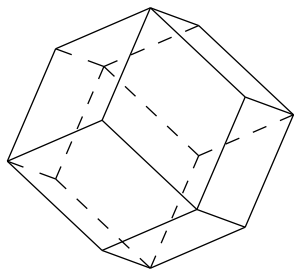

\begin{center}
\pspicture(-2,-2)(2,2)
\pnode( 0.642, 1.496){A1} \pnode( 0.642,-0.236){A2} \pnode( 1.260, 0.545){A3}
\pnode( 1.260,-1.187){A4} \pnode(-1.260, 1.187){A5} \pnode(-1.260,-0.545){A6}
\pnode(-0.642, 0.236){A7} \pnode(-0.642,-1.496){A8}
\pnode( 1.902, 0.309){B1} \pnode(-1.902,-0.309){B2} \pnode(-0.618, 0.951){B3}
\pnode( 0.618,-0.951){B4} \pnode( 0.000, 1.732){B5} \pnode( 0.000,-1.732){B6}

\psset{linecolor=black}
\ncline{B1}{A1} \ncline[linestyle=dashed]{B1}{A2} \ncline{B1}{A3}
\ncline{B1}{A4}

\ncline{B2}{A5} \ncline[linestyle=dashed]{B2}{A6} \ncline{B2}{A7}
\ncline{B2}{A8}

\ncline[linestyle=dashed]{B3}{A1} \ncline[linestyle=dashed]{B3}{A2}
\ncline[linestyle=dashed]{B3}{A5} \ncline[linestyle=dashed]{B3}{A6}

\ncline{B4}{A3} \ncline{B4}{A4} \ncline{B4}{A7} \ncline{B4}{A8}
\ncline{B5}{A1} \ncline{B5}{A3} \ncline{B5}{A5} \ncline{B5}{A7}

\ncline[linestyle=dashed]{B6}{A2} \ncline{B6}{A4}
\ncline[linestyle=dashed]{B6}{A6} \ncline{B6}{A8}

\endpspicture
\end{center}
\caption{The convex hull of $D_3$, a rhombic dodecahedron}
\label{f:rhombic}
\end{figure}

The conjecture is motivated by the fact that $D_n$ has $n(n+1)$ sides, the
largest number possible for a polytope that has the circumscribing property
\cite{Makeev:polygons}. \fig{f:rhombic} illustrates $D_3$, a standard rhombic
dodecahedron.

In this note, we will prove that every constant width body in $\R^3$ is
circumscribed by an odd number of congruent copies of $D_3$ (in a homological
sense), as is also the case in two dimensions.  In particular, we prove
Conjecture~\ref{c:main} for $n=3$, a special case which was conjectured in 1974
by Chakerian.  We also prove the following partial generalization:

\begin{theorem} \label{th:main} Every strictly convex body in $\R^3$ is
inscribed in a polyhedron which is affinely equivalent to the standard rhombic
dodecahedron.
\end{theorem}

It's not clear if the strict convexity condition is necessary.

In fact, Conjecture~\ref{c:main} and Theorem~\ref{th:main} can be generalized
further:  We can replace $D_3$ by the polyhedron
$$\begin{aligned}
P =\ & \{(x,y,z): |x| \le 1, |y| \le 1, \\
 & a|x| + a|y| + b|z| \le \sqrt{2a^2+b^2}\}
\end{aligned}$$
See Section~\ref{s:odds}.

All of these results are analogous to old results in two dimensions:  Every convex
body is circumscribed by an affinely regular hexagon and there are
homologically an odd number of them \cite{Grunbaum:hexagons}.  Instead of a
regular hexagon, we can take any centrally symmetric hexagon that circumscribes
the unit circle.

Unfortunately, for $n\ge 4$, there are homologically zero circumscribing
copies of $D_n$.  However, this does not disprove Conjecture~\ref{c:main}.

\acknowledgements

The author has learned that Makeev \cite{Makeev:rhombo} and independently
Hausel, Makai, and Szucs \cite{HMS:inscribing} have obtained similar results.

The author would like to thank Don Chakerian for suggesting the problem and
for extensive discussions, and also Bill Thurston and the
referee for useful comments.

\section{Support functions}

We establish an equivalence between constant-width bodies and antisymmetric
functions on the sphere.

Let $K$ be a convex body in $\R^n$ containing $0$, the origin. For each unit
vector $v$, let
$$f(v) = d(H_v,0),$$
where $H_v$ is the hyperplane which supports
$K$, which is orthogonal to $v$, and which is on the same side of the origin as
$v$.  The function $f$ is called the support function of $K$. The
function
$$g(v) = f(v) - 1$$
be the adjusted support function of $K$.

Conversely, if $g$ is any continuous function on the sphere $S^{n-1} \subset
\R^n$ which is strictly less than 1, and if the spherical graph of
$$f(v) = 1/(g(v)+1)$$
is convex, then $g$ is the adjusted support function of some convex body $K$,
namely the polar body of the graph of $f$.  We will call such a function $g$
pre-convex.  Moreover, $g$ is antisymmetric if and only if $K$ has constant
width 2. In conclusion, convex bodies in $\R^n$ correspond to pre-convex
functions on $S^{n-1}$ and those that have constant width 2 correspond to
antisymmetric pre-convex functions.

\begin{proposition}
Let $P$ be a polytope that circumscribes the sphere $S^{n-1}$ and
let $T$ be the set of points at which it is tangent.
Every convex body $K$ (of constant width 2) is
circumscribed by an isometric image of $P$ if and only if every continuous
(antisymmetric) function $g$ agrees with a linear function on some isometric
image of $T$. \label{p:convert}
\end{proposition}
\begin{proof}
Let $K$ be such a body and let $g$ be its adjusted support function.  The
polytope $P$ circumscribes $K$ is equivalent to the statement that $g$ vanishes
identically on $T$.  Translating $K$ is equivalent to adding a linear function
to $g$. This establishes the ``if'' direction of the proposition. It also
establishes part of the ``only if'' direction, namely for pre-convex $g$ rather
than for arbitrary continuous $g$.

Consider the set $X$ of all continuous $g$ which agree with a linear function
on some isometric image of $T$.  This set is closed under multiplication by a
scalar, and it is also a closed subset of the space of continuous functions on
$S^{n-1}$ taken with the Hausdorff topology.  If $X$ contains all pre-convex
functions, then it must be the entire space of continuous functions, because
every continuous function lies in the closure of the pre-convex functions in
this double sense.  (Any smooth function becomes pre-convex if multiplied by a
sufficiently small constant and any continuous function can be approximated by
smooth functions.) This completes the argument for the ``only if'' direction.

Both arguments also hold in the antisymmetric case.
\end{proof}

Proposition~\ref{p:convert} demonstrates that the circumscription problem for
constant-width bodies belongs to a family of questions that includes the
Knaster problem. This problem asks which finite families of points $T$ on the
unit sphere $S^{d-1} \subset \R^d$ have the property that any continuous
function from the sphere to $\R^n$ is constant on an isometric image of $T$. The
more general problem goes as follows:  Given a finite set of points $T$ on
$S^{d-1}$ and given a linear subspace $L$ of the vector space of functions from $T$
to $\R^n$, does every continuous function
$$f:S^{d-1} \to \R^n$$
admit an isometry $R$ such that $f \circ R$ lies in $L$ after restriction to
$T$? Even more generally, given any subspace $V$ of finite codimension in the
space of continuous functions on the sphere, does every continuous $f$ admit an
isometry $R$ such that $f \circ R \in V$? Of course the answer in general
depends on $V$ as well as $d$ and $n$.

If the polytope $P$ is the dual difference body $D_n$, then $T$ is the
set of vertices of the difference body of a regular simplex, also
known as the root system $A_n$.  In this case, Conjecture~\ref{c:main} is
equivalent to the assertion that for every continuous, antisymmetric $f$
on $S^{n-1} \subset \R^n$, there is a position of the root
system $A_n$ such that the restriction of $f$ is linear.

\section{Two dimensions}

The root system $A_2$ consists of six equally spaced points on the unit circle
Let $C$ be the space of all isometric images $T$ of $A_2$.  The set $C$ is a
topological circle.  It has a natural 3-dimensional vector bundle $F$ whose
fiber at each $S \in C$ is the vector space of antisymmetric functions on $T$. 
If we divide this fiber by the linear functions on $T$, the result is a new
vector bundle $E$ on $C$.  It is easy to check that the bundle $E$ is a
M\"obius bundle.

\begin{figure}[htb]
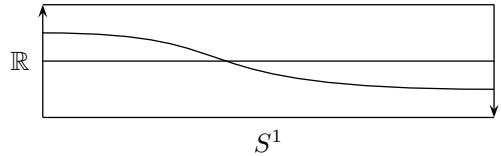

\begin{center}
\pspicture(-3,-1.1)(3,.9)
\pnode(-3, .75){A} \pnode( 3, .75){B} \pnode( 3,-.75){C} \pnode(-3,-.75){D}
\ncline{->}{D}{A} \Aput{$\R$} \ncline{A}{B} \ncline{->}{B}{C}
\ncline{C}{D} \Aput{$S^1$}
\psline(-3,0)(3,0) \psbezier(-3,.375)(0,.375)(-1.5,-.375)(3,-.375)
\endpspicture
\end{center}
\caption{A section of the M\"obius bundle}
\label{f:mobius}
\end{figure}

If $g$ is an antisymmetric, continuous function on the unit circle, it yields a
section of $F$ given by restricting $g$ to each sextuplet $T$.  In turn, one
gets a section $s$ of the bundle $E$.  We wish to know whether the section $s$
must have a zero.  Since $E$ is a M\"obius bundle, this is true (\fig{f:mobius}).

Thus we have proved that any constant-width body in the plane is 
circumscribed by a regular hexagon.  The proof is actually just the
traditional proof with some unconventional terminology.  This terminology will
be useful in the higher-dimensional cases.

\section{Three dimensions}
\label{s:three}

\begin{figure}[htb]
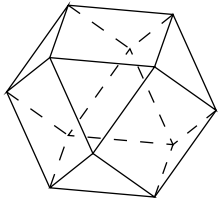

\begin{center}
\pspicture(-1.5,-1.5)(1.5,1.5)
\pnode( 0.246, 0.696){A1} \pnode( 0.819, 1.153){B1} \pnode(-0.574, 1.276){C1}
\pnode( 1.393,-0.123){A2} \pnode(-0.819, 0.579){B2} \pnode(-0.574,-0.456){C2}
\pnode(-1.393, 0.123){A3} \pnode( 0.819,-0.579){B3} \pnode( 0.574, 0.456){C3}
\pnode(-0.246,-0.696){A4} \pnode(-0.819,-1.153){B4} \pnode( 0.574,-1.276){C4}

\psline[linestyle=dashed](C1)(A1)(B1) \psline(B1)(C1)
\pspolygon[linestyle=dashed](A1)(B3)(C2)
\pspolygon(A2)(B1)(C3)
\psline[linestyle=dashed](A2)(B3)(C4) \psline(A2)(C4)
\pspolygon(A3)(B2)(C1)
\psline[linestyle=dashed](B4)(C2)(A3) \psline(A3)(B4)
\pspolygon(A4)(B2)(C3)
\pspolygon(A4)(B4)(C4)
\endpspicture
\end{center}
\caption{$A_3$ forms a cuboctahedron}
\label{f:cuboct}
\end{figure}

We wish to show that every continuous, antisymmetric function on the 2-sphere
agrees with a linear function on some isometric image of the root system $A_3$,
the vertices of a standard cuboctahedron (\fig{f:cuboct}).  The set of such
isometric images is a 3-manifold
$$M = \SO(3)/\Gamma,$$
where $\Gamma$ is the rotation group of $A_3$ (acting by right multiplication
on $\SO(3)$).  The group $\Gamma$ is the rotation group of the cube and is
isomorphic to $S_4$, the symmetric group on four letters. The manifold $M$ has
a 6-dimensional bundle $F$ which at each point is the vector space of
antisymmetric functions on the corresponding image of $A_3$.  We quotient $F$
by the linear functions to obtain a 3-dimensional bundle $E$.

We first rephrase the topological argument of the previous section in terms of
characteristic classes of vector bundles \cite{MS:characteristic}. An
$n$-dimensional bundle $B$ on an arbitrary topological space $X$ (at least a
reasonable one such as a CW complex) defines a characteristic cohomology class
$\chi(B)$ called the Euler class.  If $X$ is a closed manifold, this class is
dual to the homology class represented by the zero locus of  a generic section
of $B$.  If $B$ is orientable, then $\chi(B)$ is an element of the ordinary
cohomology $H^n(X)$. But in general
$$\chi(B) \in H^n(X,\det(B)).$$
\Ie, the Euler class lies in the cohomology of $X$ in a twisted coefficient
system, the determinant bundle of $B$.  In our case, $E$ is a non-orientable
3-plane bundle on the closed, orientable 3-manifold $M$.  Therefore
$$\chi(E) \in H^3(M,\det(E)) \cong \Z/2.$$
In other words, the Euler class $\chi(E)$ is either $0$ or $1$, depending on
whether a generic section has an even or odd number of zeroes.

\begin{theorem} The bundle $E$ has a non-trivial Euler class:
$$\chi(E) = 1 \in \Z/2.$$
\end{theorem}
\begin{proof} There are two ways to argue this.  The first way is by direct
geometric construction.  Consider the function $xyz$ on $S^2$.  It produces a
section $s$ of $E$.  The symmetry group of $xyz$, including antisymmetries, is
the same group $\Gamma$; thus, the section $s$ has the same symmetries. The
group $\Gamma$ acts on the manifold $M$ by means of symmetries that preserve or
negate $xyz$ but move some isometric image of $A_3$.  This is the {\em left}
action of $\Gamma$ on the coset space $M = \SO(3)/\Gamma$; the quotient is the
double coset space $\Gamma\backslash\SO(3)/\Gamma$.  The action has one fixed point
(coming from the identity in $\SO(3)$) and one orbit of size 3 (coming from a
rotation by 45 degrees in $\SO(3)$).  All other orbits have even order. 
An elementary calculation shows that the fixed point is a transverse
zero of the section $s$, while $s$ is non-zero on the orbit of order
3.  Thus the odd orbits make an odd contribution to the intersection
between $s$ and the zero section.  The remaining zeroes of $s$, if there
are any, lie on even-sized orbits and make an even contribution.
Thus the Euler class of $E$ is $1$ and not $0$.

The second way is by means of algebraic topology.  
Suppose that a vector bundle $V$ on a space $X$ lifts to
a {\em trivialized} bundle $\tilde{V}$ on some covering space $\tilde{X}$.
Then $V$ together with the choice of $\tilde{V}$ is called a flat
bundle.
Both $F$ and $E$ are trivial if lifted to $\SO(3)$, as well as flat
on $M$, by construction.  In general a flat bundle on a space $X$ is described
by some linear representation of the group of deck translations of $\tilde{X}$
over $X$, assuming for simplicity that the covering is regular.  In this case,
the representation $R$ of $\Gamma$  that encodes $E$ is simply the action
of $\Gamma$ on antisymmetric functions (modulo linear functions) on one copy of
the $A_3$ root system. By writing down the character of this representation, or
by writing down the representation explicitly, we can see that it is isomorphic
to $V \tensor L$, where $V$ is the  3-dimensional representation of $\SO(3)$
restricted to $\Gamma$ and $L$ is the 1-dimensional representation of $\Gamma$
coming from the sign homomorphism from $\Gamma = S_4$ to $\{\pm1\}$.
We can express this in terms of bundles with the equation
$$E \cong E_V \tensor E_L,$$
where $E_V$ and $E_L$ are
the bundles defined by the representations $V$ and $L$.  

If a flat bundle $X$ on a coset space $G/H$ is given by a representation of $H$
that is induced from $G$, it is a trivial bundle.  For example, the bundle
$E_V$ is trivial for this reason. Thus the bundle $E$ is actually three copies
of the line bundle $E_L$. It is a general property of Euler classes that if $X$
and $Y$ are two bundles, the Euler class of the direct sum is the cup product
of the Euler classes:
$$\chi(X \oplus Y) = \chi(X) \cup \chi(Y).$$
In this case we begin with the simpler Euler class
$$c = \chi(E_L) \in H^1(M,\Z/2)$$
from which we compute
$$\chi(E) = c \cup c \cup c.$$
We abbreviate $H^i(M,\Z/2)$ as just $H^i$. The cohomology group $H^1$ can be
understood as the set of homomorphisms from $\pi_1(M)$ to $\Z/2$.  In this
case, all homomorphisms factor through $\Gamma$ and $H^1 \cong \Z/2$. By this
interpretation $c$ is the same homomorphism as the one defining $L$, \ie, the
non-trivial one.  By Poincar\'e duality, $H^2 \cong \Z/2$ as well, while $H^3
\cong \Z/2$ automatically because $M$ is a closed 3-manifold. The cup product
$$\cup:H^1 \times H^2 \to H^3$$
is a non-degenerate pairing.  To determine $\chi(E)$, the only
question is whether $c \cup c$ is non-zero.  In general, if $X$ is a reasonable
topological space and $x \in H^1(X,\Z/2)$ corresponds to a homomorhism from
$\pi_1(X)$ to $\Z/2$, then $x \cup x$ vanishes if and only if the homomorphism
lifts to $\Z/4$.  One can check that the sign homomorphism of $\Gamma$ does not
lift, so $c \cup c$ is non-zero.  Therefore the Euler class
$\chi(E)$ does not vanish, as desired.
\end{proof}

Let $\Gamma_2$ be the Sylow 2-subgroup of $\Gamma$ and let
$$M_2 = \SO(3)/\Gamma_2$$
be the corresponding covering space of $M$.  Since the
covering $M_2 \to M$ has odd degree, the lift of the bundle $E$ to $M_2$ also
has odd Euler class.  This means that the theorem that every constant-width
body is circumscribed by a $D_3$ generalizes to other polyhedra $P$ with
symmetry group $\Gamma_2$, provided that the corresponding bundle $E_P$ on
$M_2$ is isomorphic to $E$, or that the corresponding representation is still
$R$.  For example, $P$ can be any of the dodecahedra mentioned in the
introduction.

\section{The bad news}
\label{s:bad}

In any dimension $n$, there is a rotation group $\Gamma$ which
preserves the $A_n$ root system and there is a manifold
$$M = \SO(n)/\Gamma$$
of positions of the root system.  The set of antisymmetric functions
modulo linear functions is a flat bundle $E$ on $M$ whose dimension
agrees with $M$.  Let
$$d = n(n-1)/2$$
be the dimension of $M$.

If $n$ is 0 or 1 modulo 4, the bundle $E$ is orientable, and its Euler class is
therefore an element of $H^d(M,\Z)$, \ie, an integer, if an orientation is
chosen.  In general, the rational Euler class of a bundle $X$ has a  Chern-Weil
formula, an expression in terms of the curvature of $X$.  Since our bundle $E$
is flat, this integral expression vanishes.  The Euler class is therefore 0.
Another way to argue this is that, as in 3 dimensions, $E$ is a sum of line
bundles.  Negating one of the line bundles yields an orientation-reversing
automorphism of $E$.   The existence of such an automorphism tells us that the
Euler class is its own negative.

If $n$ is $2$ or $3$ modulo 4, the Euler class is an element of
$$H^d(M,\Z/2) \cong \Z/2.$$
We argue that for $n \ge 4$, this number also vanishes.

\begin{proposition} For $n \ge 4$, $M$ admits a fixed-point free involution
$\sigma$ that extends to $E$.
\end{proposition}

If we accept this proposition, we are done, since whatever $\chi(E)$ is on
$M/\sigma$, it is an even multiple of it on $M$ itself.  It therefore vanishes
modulo 2.

\begin{proof}(Sketch) It suffices to find an involution $g$ in $\SO(n)$ that 
centralizes $\Gamma$ but is not in $\Gamma$.  For then the group  $\Gamma'$
generated by $\Gamma$ and $g$ would be a Cartesian product $\Gamma \times
\Z/2$, the linear representation $R$ would extend from $\Gamma$ to $\Gamma'$,
and the bundle $E$ would descend from $M$ to
$$M/g = \SO(n)/\Gamma'.$$

The group of all isometries of a simplex in $\R^n$ is the  permutation group
$S_{n+1}$.  Adding central inversion, the full isometry group of $D_n$ is
$$S_{n+1} \times \Z/2 \subset \O(n).$$
The group $\Gamma$ is an index 2 subgroup of this isometry group. The
embedding of $S_{n+1}$ in $\O(n)$ is a linear representation which
is almost the linear extension of the permutation representation
on $n+1$ letters; the difference is that a trivial summand has
been deleted.  Let $S_{n+1;2}$ be the Sylow 2-subgroup of $S_{n+1}$.
The action of $S_{n+1;2}$ on $\R^n$ can be analyzed with arcane but
standard computations.  The property of this action that we need
is that for $n \ge 4$, there are more representation endomorphisms
in $\O(n)$ (meaning isometries that commute with the action of
$S_{n+1;2}$) than those provided by the center of $S_{n+1;2} \times
\Z/2$ \cite{Guralnick:private}.  These extra endomorphisms include
orientation-preserving involutions.  The element $g$ above can be any
such involution.
\end{proof}

The author also considered the natural conjecture that every constant-width
body $K$ in $\R^4$ is circumscribed by a regular cross polytope $C$
(generalized octahedron).  Since it has four fewer sides than the polytope
$D_4$, a 2-parameter family of copies of $C$ circumscribes $K$ if $K$ is
chosen generically.  Unfortunately, another calculation shows that the set of
such circumscribing polytopes is null-homologous in $SO(4)/\Gamma$, where
$\Gamma$ is the rotation group of $C$.

Finally, a constant-width body in $\R^3$ is inscribed in homologically zero
regular dodecahedra.  Chakerian has also asked whether there is always
such a dodecahedron.

\section{Affine circumscription}
\label{s:affine}

Interestingly, the affine case of theorem~\ref{th:main} is a corollary of the
constant-width case.  For simplicity we begin with the argument in two
dimensions.  It is again an Euler class argument, except  it is more
complicated because the base space of the bundle is not compact.  In this case
a section of the bundle has a well-defined Euler class if it is proper, in the
same sense that a map between non-compact spaces may be proper.

Let $K$ be a convex body in the plane and let $H$ be a regular hexagon whose
inscribed circle has radius 1.  Let
$$G = \mathrm{GL^+}(2,\R) \ltimes \R^2$$
be the space of orientation-preserving affine transformations of the plane.
There is a map $\Phi$ from $G$  to $\R^6$ defined as follows:  For a given
affinity $\alpha$, the coordinates of $\Phi(\alpha)$ are the distances from the
lines containing the sides of $H$ to $\alpha(K)$.  If $\alpha(K)$ is on the
same side of such a line as $H$ is, the distance is taken to be negative,
otherwise it is positive. Apparently $\Phi$ is continuous.

We wish to show that for a sufficiently small open neighborhood $U$ of 0,
$\Phi^{-1}(U)$ is bounded (contained in a compact subset of $G$), for then
$\Phi$ has a well-defined degree.  More precisely, we identify $\R^6 - U$ to a
point to make the target of $\Phi$ a ball, and we extend the domain to the
one-point compactification of $G$.  Then the degree of $\Phi$ is the degree of
this modified map.

\begin{lemma} For a suitable $U$ (independent of $K$) containing the origin,
the region $\Phi^{-1}(U)$ is contained in a compact set. \label{l:compact}
\end{lemma}
\begin{proof} (Informally)
We argue that if $\alpha \in G$ is sufficiently close to infinity,
$\Phi(\alpha)$ is bounded away from $0$.  (Sufficiently close to infinity means
sufficiently near the compactification point in the one-point compactification
of $G$, or outside of a sufficiently large compact subset of $G$.) In general
an element $\alpha$ may be close to infinity if the corresponding affine image
$\alpha(K)$ has one of four properties:  It may be translated far from $H$, it
may be tiny, it may be enormous, or it may be highly anisotropic
(needle-like).  In the first three cases $\Phi(\alpha)$ is clearly bounded away
from $0$.

\begin{figure}[htb]
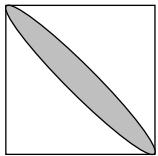

\begin{center}
\pspicture(-1,-1)(1,1)
\psframe(-1,-1)(1,1)
\rput{45}(0,0){\psellipse[fillstyle=solid,fillcolor=lightgray](0,0)(.2,1.4)}
\endpspicture
\end{center}
\caption{A needle-like ellipse inscribed in a square}
\label{f:needle}
\end{figure}

The last case is more subtle, particularly since the conclusion would not hold
if $H$ were a square rather than a hexagon (\fig{f:needle}). However, the
smallest convex body inscribed in a regular hexagon is an equilateral triangle
meeting three vertices.  This follows from the more general fact that the
smallest convex body inscribed in an arbitrary convex polygon is the convex
hull of some of the vertices.  (Such a body must touch each side and one of the
endpoints of each side is always better than points in the middle.) If
$\alpha(K)$ is so needle-like that its area is half of that of this triangle,
then $\Phi(\alpha)$ is again bounded away from $0$.
\end{proof}

Since the set $U$ in Lemma~\ref{l:compact} is independent of $K$, and since $K$
can be varied continuously, the degree of $\Phi$ is independent of $K$ as well.
Unfortunately it vanishes. However, the rotation group $\Gamma$ of $H$ acts on
$G$ and on $\R^6$, and $\Phi$ is equivariant with respect to this action. Thus
$\Phi$ represents a section of a bundle $F$ on $W = G/\Gamma$ that also
satisfies Lemma~\ref{l:compact}.

The section $\Phi:W \to F$ has an Euler class rather than a degree. To compute
it we take $K$ to be the unit circle. The zero locus of $\Phi$
is then $M = \SO(2)/\Gamma$, the manifold that appears in the constant-width
case.   Moreover, $\Phi$ is transverse to the zero section of $F$ in the
directions normal to $M$.  These directions are characterized by affinities
whose matrices are symmetric, \ie, by stretching or squeezing $K$ along
orthogonal axes.   The derivative of such a motion is radially a homogeneous
quadratic function on the boundary of the circle $K$.  The key fact to check is
that a homogeneous quadratic function is determined by its values on $A_2$, the
tangencies of the hexagon $H$.  In other words, the derivative of $\Phi$ here
is essentially restriction to $A_2$, a linear transformation which is
nonsingular for homogeneous quadratic functions.  If we quotient $F$ on $M$ by
the image under $\Phi$ of the normal bundle $NM$ of $M$, we are left with the
bundle $E$ on $M$ considered previously. Thus the Euler class of $\Phi$ on $W$
equals the Euler class of $E$ on $M$, namely $1 \in \Z/2$.

\begin{figure}[htb]
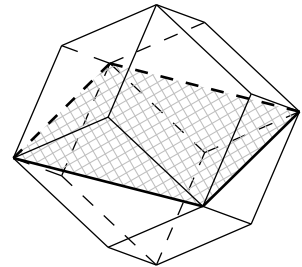

\begin{center}
\pspicture(-2,-2)(2,2)
\pnode( 0.642, 1.496){A1} \pnode( 0.642,-0.236){A2} \pnode( 1.260, 0.545){A3}
\pnode( 1.260,-1.187){A4} \pnode(-1.260, 1.187){A5} \pnode(-1.260,-0.545){A6}
\pnode(-0.642, 0.236){A7} \pnode(-0.642,-1.496){A8}
\pnode( 1.902, 0.309){B1} \pnode(-1.902,-0.309){B2} \pnode(-0.618, 0.951){B3}
\pnode( 0.618,-0.951){B4} \pnode( 0.000, 1.732){B5} \pnode( 0.000,-1.732){B6}

\psset{linecolor=black}
\ncline[linestyle=dashed]{B1}{A2}
\ncline[linestyle=dashed]{B2}{A6}
\ncline[linestyle=dashed]{B6}{A2}
\ncline[linestyle=dashed]{B6}{A6}
\ncline[linestyle=dashed]{B3}{A1}
\ncline[linestyle=dashed]{B3}{A2}
\ncline[linestyle=dashed]{B3}{A5}
\ncline[linestyle=dashed]{B3}{A6}


\pspolygon[linestyle=none,fillstyle=crosshatch,hatchwidth=.4pt,hatchsep=3pt,
    hatchangle=30,hatchcolor=lightgray](B3)(B2)(B4)(B1)

\ncline[linewidth=1pt,linestyle=dashed]{B3}{B2} \ncline[linewidth=1pt]{B2}{B4}
\ncline[linewidth=1pt,linestyle=dashed]{B1}{B3} \ncline[linewidth=1pt]{B4}{B1}

\ncline{B1}{A1}
\ncline{B1}{A3}
\ncline{B1}{A4}

\ncline{B2}{A5} 
\ncline{B2}{A7}
\ncline{B2}{A8}

\ncline{B4}{A3} \ncline{B4}{A4} \ncline{B4}{A7} \ncline{B4}{A8}
\ncline{B5}{A1} \ncline{B5}{A3} \ncline{B5}{A5} \ncline{B5}{A7}

\ncline{B6}{A4}
\ncline{B6}{A8}

\endpspicture
\end{center}
\caption{A square inscribed in a rhombic dodecahedron}
\label{f:square}
\end{figure}

This argument generalizes verbatim to three dimensions,
except that unfortunately Lemma~\ref{l:compact} no longer holds.  Among closed
convex sets inscribed in the rhombic dodecahedron $D_3$, a square, which has volume
zero, has the least volume (\fig{f:square}).  The square is the unique minimum up to
isometry.  If $K$ is strictly convex, its affine image $\alpha(K)$ is bounded
away from a square, and therefore $\Phi(\alpha)$ is again bounded away from $0$
for $\alpha$ sufficiently close to infinity.

Thus in three dimensions the Euler class of $\Phi$ is well-defined when $K$ is
strictly convex.  Moreover, a finite path $\{K_t\}_{t \in [0,1]}$ of strictly
convex bodies is strictly convex in a uniform fashion by compactness. Therefore
the Euler class of $\Phi$ does not change along such a path.  For every
strictly convex $K$ it must always equal its value when $K$ is a round
sphere, namely $1 \in \Z/2$.

\section{Odds and ends}
\label{s:odds}

Following the computations of Section~\ref{s:bad}, we did not really need the
full symmetry group of  the rhombic dodecahedron $D_3$ full symmetry group, but
only its Sylow 2-subgroup $\Gamma_2$ and the way that this subgroup permutes
its faces.  Because if we lift the bundle $E$ of Section~\ref{s:three} to an
odd-order covering of $M$, its Euler class remains non-zero. Thus the argument
applies to any other polytope which is symmetric under $\Gamma_2$, whose faces 
are permuted by $\Gamma_2$ in the same way, which is centrally symmetric, and
which circumscribes the sphere.  In particular the results hold for the
polytope $P$ described in the introduction.

It would be interesting if there were a convex body $K$ which does not affinely
inscribe in a rhombic dodecahedron.  We can obtain some information about such
a $K$ from the arguments of Section~\ref{s:affine}. It would necessarily
affinely project onto a square.  Given any sequence of strictly convex bodies
$$K_1, K_2, \ldots \to K,$$
their affine inscriptions in $D_3$ would necessarily converge to an inscribed 
square.  Affine circumscriptions of $D_3$ around each $K_n$ would converge to
an infinite parallelogram prism circumscribing $K$, and $K$ would meet all four
edges of this prism.  Otherwise some subsequence of the affine images of $D_3$
would converge to an affine image circumscribing $K$.


\end{document}